\documentclass[twoside,leqno,fleqn]{article}
\usepackage{amsmath}
\usepackage{amssymb}
\usepackage[T1]{fontenc}
\usepackage[shortlabels]{enumitem}
\usepackage[onehalfspacing]{setspace}

\newcommand{\AUTOR}{C. Ga\ss ner}
\newcommand{\TITEL}{Relationships between Principles of Choice}
\def\S{\Sigma}
\def\I{{\cal I}}
\def\G{\mathfrak{G}}

\newcommand{\bbbn}{\mathbb{N}}

\newcommand{\mbm}[1]{\mbox{\boldmath{$#1$}}}
\newcommand{\mbmss}[1]{\mbox{\scriptsize\boldmath{$#1$}}}
\newcommand{\mbmty}[1]{\mbox{\tiny\boldmath{$#1$}}}
\newcommand{\qed}{\hfill{$\Box$}} 

\newtheorem{satz}{Satz}[section] 
\newtheorem{lemma}[satz]{Lemma}
\newtheorem{proposition}[satz]{Proposition}
\newtheorem{corollary}[satz]{Corollary}
\newtheorem{theorem}[satz]{Theorem}
\newtheorem{example}[satz]{Example}
\newtheorem{defi}[satz]{Definition}

\begin{document}
\newcounter{li}

\thispagestyle{empty}
\begin{center} {\Large\bf Relationships between Principles of Choice\vspace{0.4cm}\\ in Second-Order Henkin Structures}\vspace{0.6cm}\\{\bf Christine Ga\ss ner}\footnote{I thank the organizers of the DMV Annual Meeting in Berlin in 2022. In particular, I thank Sandra Müller and Aleksandra Kwiatkowska for organizing Section S01 in Berlin.} 
 {\vspace{0.2cm}\\Universit\"at Greifswald, Germany, 2024\\ gassnerc@uni-greifswald.de}\\\end{center} 

\begin{abstract} We deal with the strength of classical second-order versions of the Axiom of Choice (AC) in second-order predicate logic (PLII) with Henkin interpretation (HPL). We use the known relationships between the so-called Zermelo-Asser axioms and the so-called Russell-Asser axioms and prove relationships between the so-called Ackermann axioms and the Zermelo-Asser axioms and between the so-called Asser axioms and the Zermelo-Asser axioms. In particular, we give the technical details of the proofs of our results presented at the DMV Annual Meeting 2022. 
\end{abstract}

\section{Introduction: Classical second-order versions of AC in HPL}\label{SectionEinf}

Motivated by \cite{1c}, the dissertation \cite{Gass84} investigates the relative strength of second-order formulas whose first-order counterparts are known as equivalent to the Axiom of Choice (denoted by AC) in ZF and other modifications of AC in {\rm HPL} (for an overview see \cite{Gass94}). Here, we continue this approach. We consider the following principles of choice discussed at the DMV Annual Meeting 2022 and the CL 2022 (cf.\,\,\cite{Gass22A} and \cite{Gass22B}) and use the definitions given in \cite{1c}, \cite{Gass94}, and \cite{Gass24A}. In formalizing these principles we will also use the following variables. For given $n\geq 1$ and $m\geq 1$, $\mbm{x}$ stands for $n$ variables $x_1,\ldots,x_n$, $\mbm{y}$ stands for $m$ variables $y_1,\ldots,y_m$, and so on. Let $A$ be a variable of sort $n$, let $C, C_1, C_2$, and $D$ be variables of sort $m$, and let $R$ and $S$ be of sort $n+m$. The set ${\cal L}^{(2)}_{\mbmss{x},D}$ contains all second-order formulas in which the individual variables in $\mbm{x}$ and $D$ occur only free, ${\cal L}^{(2)}$ contains all second-order formulas, and so on. The  following axioms  will be given  for all $n,m\geq 1$. For more details, see \cite{Gass24A}.

\vspace{0.3cm}
\noindent{\bf The set $choice^{(2)}$ of all Ackermann axioms in PL\,II} 

\nopagebreak 
\noindent\fbox{\parbox{11.7cm}{

\vspace*{0.1cm}
$\begin{array}{ll}choice^{(2)}&=_{\rm df}\bigcup_{n,m\geq 1 }choice^{n,m}\\
choice^{n,m}&=_{\rm df}\{choice^{n,m}(H)\mid H \mbox{ is a formula in ${\cal L}^{(2)}_{\mbmss{x},D}$} \}\\
choice^{n,m}(H)&=_{\rm df} 
 \forall \mbm{x} \exists D\, H(\mbm{x},D) \to \exists S \forall \mbm{x }\, H(\mbm{x}, \lambda\mbm{y}.S\mbm{x}\mbm{y})
\end{array}$
} } 

\vspace{0.3cm}
\noindent{\bf The set $choice_h^{(2)}$ of all Ackermann axioms in HPL} 

\nopagebreak 
\noindent\fbox{\parbox{11.7cm}{

\vspace*{0.1cm}
$\begin{array}{ll}choice_h^{(2)}&=_{\rm df}\bigcup_{n,m\geq 1 }choice_h^{n,m}\\
choice_h^{n,m}&=_{\rm df}\{choice_h^{n,m}(H)\mid H \mbox{ is a formula in ${\cal L}^{(2)}_{\mbmss{x},D}$} \}\\
choice_h^{n,m}(H)\!&=_{\rm df} \!
 \forall \mbm{x} \exists D H(\mbm{x},D) \!\to \exists S \forall \mbm{x } \exists D (\forall\mbm{y}(D \mbm{y} \!\leftrightarrow S\mbm{x}\mbm{y}) \land H(\mbm{x}, D )) 
\end{array}$
} } 

\vspace*{0.3cm}

\noindent{\bf The Hilbert-Ackermann axioms $AC^{n,m}(H)$ of second order} 

\nopagebreak 
\noindent\fbox{\parbox{11.8cm}{

\vspace*{0.1cm}
$AC^{n,m}(H)\!=_{\rm df}\!\forall A \exists S
(\forall \mbm{x}
 ( A\mbm{x}\! \leftrightarrow \! \exists \mbm{y}
 H(\mbm{x},\!\mbm{y})) \!\to \forall \mbm{x}
 (A\mbm{x}\!\to\! \exists!! \mbm{y}
 (H(\mbm{x},\!\mbm{y})
 \land S\mbm{x}\mbm{y})))$
}} 

\vspace{0.3cm}

\noindent{\bf The Zermelo-Asser axioms $AC^{n,m}$ of second order} 

\nopagebreak 
\noindent\fbox{\parbox{11.8cm}{

\vspace{0.1cm}
$AC^{n,m}=_{\rm df}\forall A \forall R \exists S
(\forall \mbm{x}
 ( A\mbm{x} \leftrightarrow \exists \mbm{y}
 R\mbm{x}\mbm{y}) \to \forall \mbm{x}
 (A\mbm{x}\to \exists!! \mbm{y}
 (R\mbm{x}\mbm{y}
 \land S\mbm{x}\mbm{y})))$
 }}

\vspace{0.3cm}

\noindent{\bf The Russell-Asser axioms $AC_*^{n,m}$ of second order} 

\nopagebreak
\noindent\fbox{\parbox{11.8cm}{

\vspace{0.1cm} 
$AC_*^{n,m}=_{\rm df}\forall A \forall R \exists S
(\forall \mbm{x}
 ( A\mbm{x} \! \leftrightarrow \! \exists \mbm{y}
 R\mbm{x}\mbm{y} ) $

\hfill $
 \land \forall \mbm{x}_{1} \forall \mbm{x}_{2}
 ( A\mbm{x}_1 \land A\mbm{x}_2 \land \mbm{x}_{1}\! \neq \!\mbm{x}_{2}
 \to \lnot \exists \mbm{y}
 (R\mbm{x}_{1}\mbm{y}
 \land R\mbm{x}_{2}\mbm{y}))$
 
 \hspace{3cm}$
 \to \forall \mbm{x}(A\mbm{x}
 \to \exists !!\mbm{y}(R\mbm{x}\mbm{y}
 \land S\mbm{x}\mbm{y})))$
}}

\vspace{0.3cm}

\noindent{\bf The set $choice_*^{(2)}$ of all Asser axioms of in HPL}
\nopagebreak

\noindent\fbox{\parbox{11.8cm}{

\vspace{0.1cm} 
$\begin{array}{lll}
choice_*^{(2)}&=_{\rm df}\bigcup_{m\geq 1 }choice_*^{m}\\

choice_*^{m}&=_{\rm df}\{choice_*^{m}(H)\mid H \mbox{ is a formula in ${\cal L}_C^{(2)}$} \}\\

choice_*^{m}(H)\!&=_{\rm df}\forall C(H(C)\to\exists \mbm{y}C\mbm{y}) \\

&\,\,\, \land\,\, \forall C_1\forall C_2 (H(C_1)\land H(C_2)\land C_1\neq C_2\,
\!\to \neg \exists \mbm{y}(C_1\mbm{y}\land C_2\mbm{y}))\\

&\quad \to\,\,\exists D\forall C(H(C)\to \exists !! \mbm{y}(C\mbm{y}\land D\mbm{y}))
\end{array}$
}} 

\vspace{0.4cm} 

A second-order formula $H\in {\cal L} ^{(2)}$ is ({\em generally}) {\em valid in {\rm HPL}} if every Henkin structure of second order is a model of $H$. To express this fact we use the notations $^{h}ax^{(2)}\models H$ and $\models_h H$. For ${\cal G}\subseteq {\cal L} ^{(2)}$, we write ${\cal G} \models_h H$ if any Henkin structure that is a model of each $G\in {\cal G}$ is also a model of $H$. We write $^{h}ax^{(2)}\vdash H$ or $\vdash_h H$ if $H$ is derivable from $^{h}ax^{(2)}$ with respect to the usual rules of inference. Moreover, for ${\cal G}\subseteq {\cal L} ^{(2)}$, let ${\cal G} \vdash_h H$ stands for $^{h}ax^{(2)}\cup {\cal G}\vdash H$. 

In the following, we consider second-order permutation models that we denote by $\S_0$, $\S_1$, $\S_2$, and $\S_3$. $\S_0$ is the basic Fraenkel model of second order and $\S_1$ is the ordered Mostowski-Asser model of second order. Note, that $\S_0$ is here the Fraenkel model of second order denoted by $\S_1$ in \cite{Gass94}.

For every one-sorted first-order basic structure ${\sf S}$ with the universe $I$,
the {\em predicate basic structure $\S({\sf S})$ in $ {\sf struc}_{\rm pred}^{(\rm m)}(I)$} contains all predicates $\alpha$ for which the corresponding relation $\widetilde \alpha$ belongs to the underlying structure ${\sf S}$ or the graph of a function $\mbm{g}:I^{p}\to I$ in ${\sf S}$ given by ${\rm graph}(\mbm{g})=_{\rm df}\{( \mbm{\xi}\,.\,\eta)\mid \mbm{\xi}\in I^{p} \,\,\&\,\, \mbm{g}( \mbm{\xi})= \eta\}$ is the corresponding relation $\widetilde \alpha\subseteq I^{p+1}$. Here, $(\mbm{\xi}\,.\,\eta)$ is the tuple $(\xi_1, \ldots, \xi_{p}, \eta)$. Let ${\rm auto}({\sf S})$ be the group of all automorphisms of ${\sf S}$.

\section{Relationships between these versions of choice}\label{AbschnittAC}\label{SectionTheAxiom} 
\setcounter{satz}{0}

In analogy to \cite{Russ}, Asser defined $AC_*^{1,1}$ in \cite{1c} and he conjectured that $AC_*^{1,1}$ is weaker than $AC^{1,1}$ in {\rm HPL}. This conjecture was proved in \cite{Gass84} (cf.\,\,also \cite{Gass94}). By \cite[Section IV, Satz 2.2.1 a)]{Gass84}, $\vdash_h AC^{n,m}\to AC^{1,1}$ holds and, by \cite[Section IV, Satz 2.2.3]{Gass84}, we have $\vdash_h AC^{1,1}\to AC_*^{1,2}$. The ordered Mostowski-Asser model $\S_1$ of second order is a model of $\neg AC_*^{1,2}$ (cf.\,\,\cite[Section V, Satz 12.3]{Gass84}) and consequently a model of $ \neg AC^{n,m}$. Moreover, $\S_1$ is a model of $AC_*^{n,1}$ for all $n\geq 1$ (cf.\,\,\cite[Section V, Satz 12.5]{Gass84}). 

\begin{theorem}[Gaßner {\normalfont \cite{Gass84}}]\label{RussellZermeloHPL} In \,{\rm HPL}, $AC_*^{1,1}$ is weaker than $AC^{n,m}$ for any $n,m\geq 1$.
\end{theorem}

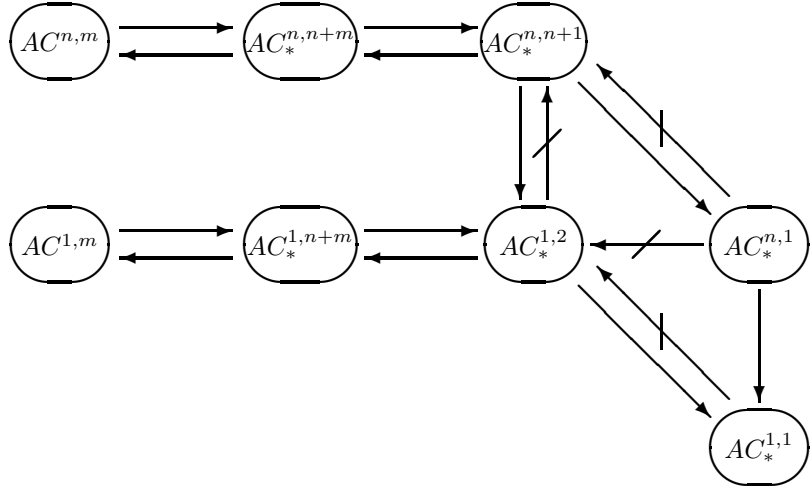
\begin{figure}[ht] \unitlength1cm

\begin{picture}(11.8,6.8) \thicklines
\put(-1,0){
\put(-3,2.7)
{
\put(5.75,3.38){\vector(1,0){1.5}}
\put(7.25,3.02){\vector(-1,0){1.5}}
\put(5.75,0.68){\vector(1,0){1.5}}
\put(7.25,0.32){\vector(-1,0){1.5}}

\put(4.95,3.2){\oval(1.3,1)\makebox(0,0){$AC^{n,m}$}}
\put(4.95,0.5){\oval(1.3,1)\makebox(0,0){$AC^{1,m}$}}
}

\put(1.75,0)
{
\put(4.25,6.08){\vector(1,0){1.5}}
\put(5.75,5.72){\vector(-1,0){1.5}}
\put(4.25,3.38){\vector(1,0){1.5}}
\put(5.75,3.02){\vector(-1,0){1.5}}
\put(6.32,5.3){\vector(0,-1){1.5}}
\put(6.68,3.8){\vector(0,1){1.5}}
\put(6.5,4.37){\line(1,1){0.36}}

\put(7.1,5.35){\vector(1,-1){1.75}}
\put(9.1,3.85){\vector(-1,1){1.75}}
\put(8.2,4.5){\line(0,1){0.48}}

\put(7.1,2.65){\vector(1,-1){1.75}}
\put(9.1,1.15){\vector(-1,1){1.75}}
\put(8.2,1.8){\line(0,1){0.48}}

\put(9.5,2.6){\vector(0,-1){1.5}}
\put(8.75,3.2){\vector(-1,0){1.5}}
\put(7.82,3.02){\line(1,1){0.36}}

\put(3.4,5.9){\oval(1.5,1)\makebox(0,0){$AC^{n,n+m}_{*}$}}
\put(3.4,3.2){\oval(1.5,1)\makebox(0,0){$AC^{1,n+m}_{*}$}}
\put(6.5,5.9){\oval(1.4,1)\makebox(0,0){$AC^{n,n+1}_{*}$}}
\put(6.5,3.2){\oval(1.3,1)\makebox(0,0){$AC^{1,2}_{*}$}}
\put(9.5,3.2){\oval(1.3,1)\makebox(0,0){$AC^{n,1}_{*}$}}
\put(9.5,0.5){\oval(1.3,1)\makebox(0,0){$AC^{1,1}_{*}$}}
}
}
\end{picture}
\caption{The relationships between the Zermelo-Asser axioms and the Russell-Asser axioms in {\rm HPL} ($n>1,m\geq 1$)}\label{RelZermRuss}
\end{figure}

We know that, for any $n,m,k\geq 1$, the formula $AC_*^{n,n+k}\leftrightarrow AC^{n,m}$ holds in {\rm HPL} (cf.\,\,Gaßner \cite[Section IV, Satz 2.2.3]{Gass84}) and, in particular, we have further relationships (cf.\,\,\cite{Gass94}) as illustrated in Figure \ref{RelZermRuss} where $H_{1} \longrightarrow H_{2}$ means that we have $\vdash_h H_{1}\to H_{2}$. The corresponding relationships follow from \cite[Section IV, Satz 2.2.3]{Gass84}, \cite[Section IV, Satz 2.2.1]{Gass84}, and \cite[Section IV, Satz 2.3.1]{Gass84}. 
$\mbox{$H_{1}\,\, \not\!\!\longrightarrow H_{2}$}$ means that $H_{1}$ does not imply $H_{2}$ in {\rm HPL} and that, in the case that $ H_{2}\to H_{1}$ holds in HPL, the formula $H_1$ is really weaker than $H_2$ in HPL. First, the second-order Mostowski-Asser model $\S_1$ is a model of $^hax^{(2)}\cup\{AC_*^{n,1}\,\,\land\,\,AC_*^{1,1}\,\,\land\,\, \neg AC_*^{1,2}\}$ ($n>1$) \cite[Section V]{Gass84}. Second, there is a model of $^hax^{(2)}\cup\{AC^{1,1}\,\,\land\,\, \neg AC^{2,1}\}$. It is the second model in \cite{Gass84}. For more details, see \cite[Section V, Satz 2.5]{Gass84}, \cite[Section V, Satz 2.9]{Gass84}, and \cite[Section IV, Satz 2.2.2 d)]{Gass84}.
Third, it can also be shown that $AC^{n,m}$ and thus all statements considered in Figure \ref{RelZermRuss} are true in the basic Fraenkel model $\S_0$. This means that $AC^{n,m}$ is independent of the weaker statements in the diagram.

\subsection{Ackermann and Asser axioms imply Zermelo-Asser axioms}

For proving $choice^{n,m}\vdash_h Ch^{(2)}_{n,m}$, Asser \cite{1c} suggested to substitute the formula $ \exists!! \mbm{y}D\mbm{y} \land \exists \mbm{y}(R \mbm{x}\mbm{y}\land D\mbm{y})$ for the placeholder $H(\mbm{x},D)$ in the formula $choice^{n,m}(H)$.
 In analogy to this, we can use
\[H(\mbm{x},D )=A \mbm{x} \to \ \exists!! \mbm{y}D\mbm{y} \land \exists \mbm{y}(R \mbm{x}\mbm{y}\land D\mbm{y})\] for obtaining the following.
\setcounter{equation}{0}
\begin{proposition}[Asser {\normalfont \cite{1c}}]\label{ChoiceToAC} For $n,m\geq 1$, we have
\[choice^{n,m}\vdash_h AC^{n,m}.\]
\end{proposition}

For proving that Asser axioms imply Russell-Asser axioms, we can replace, as mentioned in \cite{1c}, the formula $H( C)$.
Let this formula be defined as follows. \[H( C)=\,\underbrace{A\mbm{x} \land \underbrace{\forall \mbm{y}(C\mbm{y}\leftrightarrow R \mbm{x}\mbm{y})}_{F (\mbmss{x},C)}}_{G(\mbmss{x} , C)\,}\,\land \exists \mbm{y}C\mbm{y}.\] 

\setcounter{equation}{0} 

\begin{proposition}[{Asser \normalfont{\cite{1c}}}]\label{ChoiceToAC_Disj} For $n,m\geq 1$, we have
\[ choice_*^{m}\vdash_h AC_*^{n,m}.\]
\end{proposition}

\begin{example}[Validity in HPL]\label{Ein_anderer Beweis} 
The last part of the proof can be done as follows. We start with the assumption that

$\vdash_h H_0\to \exists D \forall C(A\mbm{x} \land F(\mbm{x} ,C)\to \exists !! \mbm{y}(C\mbm{y}\land D\mbm{y}))$ \hfill {\em (e)}

\noindent has been shown for $H_0= H_1\land H_2\land choice_*^{m}(H(C))$ given by 

 $\forall A \forall R \exists S
(\forall \mbm{x}
 ( \underbrace{A\mbm{x} \leftrightarrow \exists \mbm{y}
 R\mbm{x}\mbm{y}}_{H_1} ) $

 \hspace{1cm}$ \, \land \,\forall \mbm{x}_{1} \forall \mbm{x}_{2}
 ( \,\underbrace{A\mbm{x}_1 \land A\mbm{x}_2 \land \mbm{x}_{1} \neq \mbm{x}_{2}
 \to \lnot \exists \mbm{y}
 (R\mbm{x}_{1}\mbm{y}
 \land R\mbm{x}_{2}\mbm{y})}_{H_2}\,)$.

In the following, let $\S=(J_n)_{n\geq 0}$ be any Henkin-Asser structure, $\rho$ be an $(n+m)$-ary predicate in $J_{n+m}$, and $\alpha$ be an $n$-ary predicate in $J_n$ such that $\S\models_{f\langle{A\atop \alpha}{R\atop \rho} \rangle} H_0$ holds for all assignment $f$ in ${\rm assgn}(\S)$. Now, let $f$ be an assignment in ${\rm assgn}(\S)$. For each $\mbm{\xi}\in J_0^n$, let the $m$-ary predicate $\gamma_{\mbmss{\xi}} $ be given by $\gamma_{\mbmss{\xi}}=\alpha_{\S, R \mbmss{x}\mbmss{y},\mbmss{y},\langle{R\atop \rho}{\mbmty{x}\atop \mbmty{\xi}} \rangle}$. Because of $\S\models_{f\langle{\mbmty{x}\atop \mbmty{\xi}}{R\atop \rho}{C\atop\gamma_{\mbmty{\xi}}} \rangle}\forall \mbm{y} ( R\mbm{x}\mbm{y} \leftrightarrow C\mbm{y})$ 
and {\rm (e)} there is a predicate $\delta\in J_m$ with 

\vspace{0.2cm}

 $\S\models_{f\langle{A\atop \alpha} {\mbmty{x}\atop \mbmty{\xi}}{C\atop\gamma_{\mbmty{\xi}}}{D\atop\delta} \rangle}A\mbm{x}\to \exists !! \mbm{y}(C\mbm{y}\land D\mbm{y})$ \hfill{\rm (f)}

\vspace{0.2cm}

\noindent for all $\mbm{\xi}\in J_0^n$. Moreover, by definition of $\gamma_{f(\mbmss{x})}$,
we have 

\vspace{0.2cm}

$\S\models_{f\langle{R\atop \rho} {C\atop\gamma_{f(\mbmty{x})}} \rangle} R\mbm{x}\mbm{y} \leftrightarrow C\mbm{y}$ 

\vspace{0.2cm}

\noindent and, thus,

\vspace{0.2cm}

$\S\models_{f\langle{R\atop \rho} {C\atop\gamma_{f(\mbmty{x})}}{D\atop\delta} \rangle} R\mbm{x}\mbm{y} \land D\mbm{y} \leftrightarrow C\mbm{y}\land D\mbm{y}$. \hfill {\rm (g)}

\vspace{0.2cm}

\noindent Moreover, let $\sigma=\alpha_{\S, R \mbmss{x}\mbmss{y} \land D\mbmss{y},(\mbmss{x}\,.\,\mbmss{y}),\langle{R\atop \rho}{D\atop \delta} \rangle}$. Then, 

\vspace{0.2cm}

$\S\models_{f\langle{R\atop \rho}{S\atop \sigma} {D\atop\delta} \rangle} R\mbm{x}\mbm{y} \land D\mbm{y}\leftrightarrow R\mbm{x}\mbm{y} \land S\mbm{x}\mbm{y}$ \hfill {\rm (h)}

\vspace{0.2cm}

\noindent and, consequently, by {\rm (g)} and {\rm (h)} there hold

\vspace{0.2cm}

$\S\models_{f\langle{R\atop \rho}{S\atop \sigma} {C\atop\gamma_{f(\mbmty{x})}}{D\atop\delta} \rangle}C\mbm{y}\land D\mbm{y} \leftrightarrow R\mbm{x}\mbm{y} \land S\mbm{x}\mbm{y} $

\vspace{0.2cm}

\noindent and

\vspace{0.2cm}

$\S\models_{f\langle{R\atop \rho}{S\atop \sigma} {C\atop\gamma_{f(\mbmty{x})}}{D\atop\delta} \rangle} \exists !! \mbm{y}(C \mbm{y}\land D\mbm{y}) \to \exists !! \mbm{y}(R \mbm{x}\mbm{y}\land S\mbm{x}\mbm{y})$. 

\vspace{0.2cm}

\noindent Because of {\rm (f)} we have 

\vspace{0.2cm}

$\S\models_{f\langle{A\atop \alpha} {R\atop \rho}{S\atop \sigma} \rangle}A \mbm{x}\to \exists !! \mbm{y}(R \mbm{x}\mbm{y}\land S\mbm{x}\mbm{y})$.

\vspace{0.2cm}

\noindent This means that 

\vspace{0.2cm}

$ \models_h H_0\to \exists S\forall \mbm{x}(A\mbm{x} \to \exists !! \mbm{y}(R \mbm{x}\mbm{y}\land S\mbm{x}\mbm{y}))$.
\end{example}

In the same way, we get the following proposition if we use 
\[H(C)=(A\mbm{x}_0 \land \forall \mbm{x}\forall \mbm{y}( R \mbm{x}\mbm{y} \land \mbm{x} = \mbm{x}_0 \leftrightarrow C\mbm{x}\mbm{y}))\land \exists \mbm{x} \exists \mbm{y} C \mbm{x} \mbm{y}.\] 
\begin{proposition}\label{PropChoiceToAC} For $n,m\geq 1$, we have
\[ choice_*^{n+m}\vdash_h AC^{n,m}.\]
\end{proposition}

 Thus, by Theorem \ref{RussellZermeloHPL}, the Russell-Asser (1,1)-axiom $AC_*^{1,1}$ is weaker than $choice_*^{n+m}(H)$ ($n,m\geq 1$).

\begin{corollary}\label{CoroChoiceAndAC}
There is a model of $^{h}ax^{(2)} \cup \{AC_*^{1,1} \land \neg choice_*^2 \land \neg choice^{1,1}\}$.
\end{corollary}

\subsection{Zermelo-Asser (1,1)-axiom does not imply Asser axioms}\label{Section_Str_S2}

By using ZF as metatheory, we will now define a further permutation model denoted by $\S_2$ in order to show that $choice_*^1$ is {\rm HPL}-independent of $AC^{1,1} $. 

 In the following part of this section, let $N_1=\{0,1\}$ and each product $K\times N_1$ stand for $(K\times N_1)$ (e.i.\,\,we sometimes omit the outer parentheses of subterms of the form $(K\times N_1)$). Let the domain $I$ of individuals be given by $I=\bbbn\times N_1$ and let $\mbm{t}_n = \{((n,0),(n,1)),((n,1),(n,0))\}$ for each $n\in \bbbn$. Moreover, let $\mbm{t}\subseteq I^2$ be defined by $\mbm{t}=\bigcup_{n\in \bbbn}\mbm{t}_n $. 

The Henkin structure $\S_2$ will be the result of a combination of special technical features of the basic Fraenkel model $\S_0$ of second order with the features of the {\em second-order Fraenkel-Zermelo model} $\S(I,{\rm auto}(I;\emptyset;\emptyset;(\mbm{t}_n)_{n\geq 0}),\I_{\sf 0}^I)$. The Fraenkel-Zermelo model of second order is derived from a permutation model of ZFA that Thomas J.\,\,Jech called the second Fraenkel model (cf.\,\,\cite{Jech73}). It was introduced by Abraham A.\,\,Fraenkel (cf.\,\cite{FA22a,FA37}) to discuss a weak form of the axiom of choice for a countable set of finite sets (cf.\,\cite[pp.\,1--2]{FA37}). By an idea of Ernst Zermelo, this model can be built on a countable set of unordered pairs (cf.\,\cite[Footnote 1, p.\,254]{FA22a}). An extension of the following model is introduced in \cite[Section V.8]{Gass84}. The mentioned models are also used in \cite{Mackereth}.  $\mbm{t}$ is a binary relation on $I$. Let $\sf S$ be the one-sorted first-order structure $(I;\emptyset;\emptyset;\{\mbm{t}\})$. For defining $\S_2$, we will use the structure $\sf S$ as basic structure. Let $\G={\rm auto}({\sf S})$.
 
 \begin{defi}[The model $\S_2$ of second order]\label{Model2} \hfill
Let $\S_2$ be the Henkin-\linebreak Asser structure $\S(I,\G,\I_{\sf 0}^I)$.
\end{defi}
 
 Thus, $\S({\sf S})$ is a substructure of $\S_2$. $\sf S$ corresponds to --- this means that it can be considered to be --- an infinite (locally finite) directed graph composed of an infinite number of disjoint complete digraphs $(V_n,E_n)$ of order two with $n\in \bbbn$ (for definitions see, e.g., \cite[pp.\,209--281]{Diestel17}). $V_n$ is the set $\{v_n,w_n\}$ of the vertices $v_n=(n,0)$ and $w_n=(n,1)$ and $E_n$ is the set $\{(v_n,w_n),(w_n,v_n) \}$ of two directed edges. This means that $\G$ can be characterized by $\G =\{\pi \in \G_{\sf 1}^I \mid (\exists \phi_\pi \in \G_{\sf 1}^\bbbn)( \forall n \in\bbbn)(\exists \pi_n \in \G_{\sf 1}^{N_1})(\forall m\in N_1)(\pi(n, m) =(\phi_\pi (n),\pi_{n}( m))\}$. Each permutation in $ \G$ maps each of the subgraphs $(V_n,E_n)$ into a complete subgraph of order two. 
 
Let $\G_0^\bbbn$ contain only the identity function $\pi_{\rm id}\in \G_1^\bbbn$.
If we would replace $\G_1^\bbbn$ in ``$ \phi_\pi \in\G_1^\bbbn$'' by $\G_0^\bbbn$, then $\phi_\pi $ could be only $\pi_{\rm id}$ and we would get $ {\rm auto}(I;\emptyset;\emptyset;(\mbm{t}_n)_{n\geq 0})$ instead of $\G$. Clearly, each permutation in the group $ {\rm auto}(I;\emptyset;\emptyset;(\mbm{t}_n)_{n\geq 0})$ allows only to interchange the vertices within the single complete subgraphs $(V_n,E_n)$, i.e., for an arbitrary number of pairs (edges), the components of pairs of two vertices can be interchanged and each pair (vertex) in any set $V_n$ ($n\in \bbbn$) is mapped into the same set $V_n$.

This means that we replace each individual $\nu$ of the domain $\bbbn$ of $\S_0$ by two individuals $(\nu,0)$ and $(\nu,1)$ in order to get the domain $\bbbn\times N_1$ of $\S_2$. Each permutation $\pi$ in $\G$ results from transforming some permutation $\phi \in \G_1^\bbbn$ such that the extended permutation $\pi^*$ in $\G^*$ transfers, for instance, the predicate $\alpha$ given by $\widetilde \alpha=\{\nu\}\times N_1$ for any $\nu\in \bbbn$ to the predicate $\beta$ given by $\widetilde \beta=\{\phi(\nu)\}\times N_1$. For any $\alpha \in J_n(I,\G,\I^I_0)$, there is a $P\in \I^I _0$ with ${\rm sym}_\G(\alpha) \supseteq \G(P)$. Let $K=\{\nu\mid\exists \mu ((\nu,\mu)\in P)\}$. Then, $\G(P)=\G(K\times N_1)$ holds and we have ${\rm sym}_\G(\alpha) \supseteq \G(K\times N_1)$. Here, let $K^{\rm c}=\bbbn\setminus K$ for all $K\subseteq \bbbn$ and $M^{\rm c}=I\setminus M$ for all $M\subseteq I$. Consequently, for all relations $\widetilde \alpha$ in $\widetilde J_1(I,\G,\I_{\sf 0}^I)$ there is a finite set $K \subseteq \bbbn$ such that either $\widetilde \alpha \subseteq K\times N_1$ or $K^{\rm c}\times N_1 \subseteq \widetilde \alpha $. For binary relations $\widetilde \rho$ in $\widetilde J_2(I,\G,\I_{\sf 0}^I)$ and $\widetilde \delta=\{(\xi,\xi)\mid \xi\in I\}$, there is a finite set $K \subseteq \bbbn$ such that either $\widetilde \rho \cap \widetilde \delta \subseteq (K\times N_1)^2$ or $(K^{\rm c}\times N_1)^2 \cap \widetilde \delta \subseteq \widetilde \rho $. Moreover, then we also have either $\widetilde \rho \setminus \widetilde \delta \subseteq (K\times N_1)^2$ or $(K^{\rm c}\times N_1)^2 \setminus \widetilde \delta \subseteq \widetilde \rho $ and in particular the following.

\begin{lemma}\label{ToIdealOfS2} Let $\S_2=\S(I,\G,\I^I_0)$ and $\xi$ be any individual $(\nu_0, \mu_0)$ in $I$. For any binary predicate $\rho$ in $\S_2$, there is some finite $K\subseteq \bbbn$ such that ${\rm sym}_\G(\rho) \supseteq \G(K\times N_1)$ is satisfied and consequently the following equivalence holds for all $\nu\in \bbbn\setminus (K\cup \{\nu_0\})$ and all $\mu\in N_1$.
\[\mbox{$(\xi,(\nu,\mu))\in \widetilde \rho$ iff } \{\xi\} \times (( \bbbn\setminus (K\cup \{\nu_0\}))\times N_1) \subseteq \widetilde\rho.\]
\end{lemma}
\noindent{\bf Proof.} Let $\rho \in J_2(I,\G,\I^I_0)$. Then, there is a $P\in \I^I _0$ with ${\rm sym}_\G(\rho) \supseteq \G(P)$ and ${\rm sym}_\G(\rho) \supseteq \G(K\times N_1)$ for $K=\{\nu\mid\exists \mu ((\nu,\mu)\in P)\}$. Let $\nu_0\in \bbbn$, $\nu_1,\nu_2 \in \bbbn\setminus (K\cup \{\nu_0\})$, $\nu_1\not=\nu_2$, and $\mu\in N_1$. Then, there are three permutations $\pi_1,\pi_2,\pi_3 \in \G(K\times N_1)$ with (a), (b), and (c), respectively.

\vspace{0.2cm}

\begin{tabular}{clll}\\
 (a) &$\pi_1(\xi)= \xi$, & \qquad $\pi_1((\nu_1,0))=(\nu_1,1)$, & \qquad $\pi_1((\nu_1,1))=(\nu_1,0)$\\
(b) &$\pi_2(\xi)= \xi$, & \qquad $\pi_2((\nu_1,0))=(\nu_2,0)$,& \qquad $\pi_2((\nu_2,0))=(\nu_1,0)$\\
(c) &$\pi_3(\xi)= \xi$, & \qquad $\pi_3((\nu_1,0))=(\nu_2,1)$,& \qquad $\pi_3((\nu_2,1))=(\nu_1,0)$\\\\
\end{tabular}

\vspace{0.2cm}

\noindent By (a), 
$(\xi,(\nu_1,\mu))\in \rho$ holds iff $(\xi,(\nu_1,1-\mu))\in \rho$ holds. Consequently, by (b) and (c), we have

\vspace{0.2cm}

$(\xi,(\nu_1,\mu))\in \rho$ iff $(\xi,(\nu_2,0))\in \rho$ as well as $(\xi,(\nu_2,1))\in \rho$.
\hfill\qed

\vspace{0.2cm}

 Whereas Zermelo's Axiom of Choice is not true in the basic Fraenkel model of ZFA and in the second Fraenkel model of ZFA (cf.\,\,\cite{Jech73}), the Zermelo-Asser axiom $AC^{1,1}$ is true in $\S_0$ and even in $\S_2$ because we have the following by Lemma \ref{ToIdealOfS2}. 

\begin{proposition}\label{ACinS2} $\S_2 \models AC^{1,1} $.
\end{proposition}

\noindent{\bf Proof.} Let $\I$ be the ideal $\I^ I_0$ and $\S=\S_2$. Then, $\S=\S(I,\G,\I)$. Let $\rho\in J_2(I,\G,\I)$ and $\alpha\in J_1(I,\G,\I)$ satisfy $\S\models_{\langle{R\atop \rho}{A\atop\alpha}\rangle} \forall x ( Ax \leftrightarrow \exists y Rx y)$. By Lemma \ref{ToIdealOfS2}, there is a finite subset $K\subseteq \bbbn$ such that ${\rm sym}_\G(\rho)\supseteq\G(K\times N_1 )$. Now, let such a set $K$, $M= K\times N_1 $, and $M^{\rm c}= I\setminus M$ be given. The statement $(\xi, \eta)\in \rho$ implies the statement $(\pi(\xi), \pi(\eta))\in \rho$ for each $\pi \in \G(M)$. Consequences are ${\rm sym}_\G(\alpha)\supseteq\G(M)$ and, because of ${\rm sym}_\G(\chi_{M^{\rm c}})\supseteq\G(M)$, ${\rm sym}_\G(M^{\rm c}\cap \alpha)\supseteq\G(M)$. ${\rm sym}_\G(\alpha)\supseteq\G(M)$ also means that we have either $\alpha \cap M^{\rm c}=\emptyset$ or $M^{\rm c} \subseteq \alpha $. Let $\zeta_0= \min_{\rm lexico}M^{\rm c}$. Then, $\zeta_0=(\min (\bbbn\setminus K),0)$. Moreover, let $L= \{\min (\bbbn\setminus K)\}\times N_1 $ and $\zeta_1= \min_{\rm lexico}(M^c\setminus L)$. Then, $\zeta_1=(\min (\bbbn\setminus (K\cup \{\min (\bbbn\setminus K)\})),0)$. For all $\xi=(\nu,\mu)\in \alpha $, let 
\begin{eqnarray*} M_{\xi}\hspace{0.1cm}&=&\{\eta \in M\mid (\xi,\eta )\in \rho\} \setminus \{\nu\}\times N_1 , \\
 \!\!\! \!\! M^{\rm c}_{\xi}\,=_{\rm df}\,(M^{\rm c})_{\xi}\hspace{0.1cm}&=&\{\eta \in M^{\rm c}\mid (\xi,\eta )\in \rho \} \setminus \{\nu\}\times N_1 , 
\end{eqnarray*}

\noindent let
\begin{eqnarray*} \hspace*{1.3cm}
\phi(\xi)&=&\left\{\begin{array}{ll} \min_{\rm lexico}M^{\rm c}_{\xi} & \mbox{ if } M^{\rm c}_{\xi}\not= \emptyset, \\ \min_{\rm lexico}M_{\xi} &\mbox{ if } M^{\rm c}_{\xi}= \emptyset \mbox{ and } M_{\xi}\not= \emptyset,\\
 \xi& \mbox{ if } M^{\rm c}_{\xi}= \emptyset \mbox{ and } M_{\xi}= \emptyset\mbox{ and } (\xi, \xi)\in\rho,\\ (\nu,1-\mu)& \mbox{ if } M^{\rm c}_{\xi}= \emptyset \mbox{ and } M_{\xi}= \emptyset\mbox{ and } (\xi, \xi)\not\in\rho, \end{array}\right. 
\end{eqnarray*}
and let $\sigma$ be the predicate in ${\rm pred}_2(I)$ such that, for all $\xi,\eta\in I$, $\sigma (\xi,\eta)=true $ holds iff $ \phi(\xi)=\eta$ holds. Let $\phi_\sigma$ be the relation $\widetilde \sigma$ which implies $\phi_\sigma=\phi$. We want to show that $\sigma\in J_2(I,\G,\I)$ and $\S\models_{\langle{R\atop \rho}{A\atop\alpha}{S\atop\sigma}\rangle} \forall x (Ax \to \exists!! y (Rx y \land Sx y))$ hold. That means that it is enough to show that there is a finite set $P_\sigma\subseteq I$ such that ${\rm sym}_\G(\sigma)\supseteq\G(P_\sigma)$ and that $\phi_\sigma:\widetilde\alpha\to I$ is a choice function on $\alpha$.

Since the statement $(\xi, \phi_\sigma(\xi))\in \rho$ implies the statement $(\pi(\xi), \pi(\phi_\sigma(\xi)))\in \rho$ for each $\pi \in \G(M)$, we can derive the properties of a choice function for $\phi_\sigma$ as follows.

\noindent First, let $M^{\rm c}\cap \alpha\not= \emptyset$ and thus $M^{\rm c}\subseteq \alpha$. For each given $\rho$, only one of the following mutually exclusive cases --- (a), (b), (c), and (d) --- is possible regardless of which $\xi_0\in M^{\rm c} $ we choose. This means that we will see that all individuals $\xi \in M^{\rm c}$ satisfy the same case condition and that, for the fixed $\rho$, the function values of $\phi_\sigma$ are consequently determined in the same way for each $\xi \in M^{\rm c}$.

\begin{list}{(\alph{li})}{\usecounter{li}}
\item {$\phi_\sigma(\xi_0) \in M_{\xi_0}^{\rm c}$ holds for some $\xi_0\in M^{\rm c} $.}

Let $\xi_0=(\nu_0,\mu_0)$ and $\eta_0=\phi_{\sigma}(\xi_0)$. Then, $\eta_0\in (\bbbn\setminus (K\cup \{\nu_0\}))\times N_1$. This means, by Lemma \ref{ToIdealOfS2}, that $(\xi_0,\eta_0)\in \{\xi_0\} \times (M^{\rm c}\setminus \{\nu_0\}\times N_1 )\subseteq \rho$. Moreover, for each $\xi=(\nu,\mu)\in M^{\rm c} $ and any $\eta\in M^{\rm c}\setminus \{\nu\}\times N_1 $, there are a $\pi \in \G(M)$, a $ \phi_{\pi}\!\in \G^\bbbn_1$, and a $ \pi_{\nu_0}\in \G^{ N_1 }_1$ with  $\pi(\xi_0)=(\phi_{\pi}(\nu_0),\pi_{\nu_0}(\mu_0))= \xi$ and $\pi(\eta_0)=\eta$. The latter condition is possible since $\eta_0\not\in \{\nu_0\}\times N_1 $ implies $\pi(\eta_0)\not\in \{\phi_{\pi}(\nu_0)\}\times N_1 $. Thus, for $\xi\in M^{\rm c}$, we have $\{\xi\} \times (M^{\rm c}\setminus \{\nu\}\times N_1 )\subseteq \rho$ and, consequently, $M^{\rm c}_{\xi}=M^{\rm c}\setminus \{\nu\}\times N_1 $ which means that $M^{\rm c}_{\xi}\not= \emptyset$. This implies, by definition, $\xi\in\alpha$ and thus, for all $\xi\in M^{\rm c}$,
\[\phi_\sigma(\xi)=\min{\!}_{\rm lexico}M^{\rm c}_{\xi}=\left\{\begin{array}{ll}\zeta_0&\mbox{ if $\xi\not \in L$},\\ 
\zeta_1&\mbox { if $\xi\in L$}.\end{array}\right.\] 

\item {$\phi_\sigma(\xi_0) \in M_{\xi_0}$ holds for some $\xi_0\in M^{\rm c} $.}

Let $\xi_0=(\nu_0,\mu_0)$. By definition of $\phi_\sigma$, we have $M^{\rm c}_{\xi_0}=\emptyset$. Moreover, we have $M_{\xi_0}=(\{\xi_0\}\times M)\cap \rho$ since $\{\nu_0\}\times N_1\subseteq M^{\rm c}$ and there is an $\eta_0\in M_{\xi_0}$ with $(\xi_0,\eta_0)\in \rho$. Let $\xi=(\nu,\mu) \in M^{\rm c}$. By the consequences of (a), $\phi_\sigma(\xi) \in M^{\rm c}_\xi $ cannot be satisfied. Since for all $\eta\in M_{\xi_0}$ there is a $\pi \in \G(M)$ with $\pi(\xi_0)=\xi$ and $\pi (\eta)=\eta$, we have $(\xi,\eta)\in \rho$ which means that $(\{\xi\}\times M)\cap \rho\not=\emptyset$ and consequently $M_{\xi}\not=\emptyset$ hold. Since there is also a $\pi \in \G(M)$ with $\pi(\xi)=\xi_0$ and $\pi (\eta)=\eta$ for all $\eta \in M_{\xi}$, there holds even $M_{\xi}=M_{\xi_0}$ and thus, for all $\xi\in M^{\rm c} $, there holds
\[\phi_\sigma(\xi)=\min{\!}_{\rm lexico}M_{\xi_0}.\] 

\item {$\phi_\sigma(\xi_0) =\xi_0$ holds for some $\xi_0\in M^{\rm c} $.}

Then, by definition of $\phi_\sigma$, we have $(\xi_0, \xi_0)\in\rho$ and $M_{\xi_0}= M_{\xi_0}^{\rm c}=\emptyset$. Thus ${\rm sym}_\G(\rho)\supseteq \G(M)$ implies $(\xi, \xi)\in\rho$ for any $\xi=(\nu,\mu) \in M^{\rm c}$. Moreover, by the consequences of (a) and (b), $\{\xi\}\times (M_{\xi}\cup M_{\xi}^{\rm c})= (\{\xi\}\times (I\setminus \{\nu\}\times N_1 ))\cap \rho=\emptyset$ and thus, for all $\xi\in M^{\rm c} $, 
\[\phi_\sigma(\xi) =\xi.\] 

\item {$\phi_\sigma(\xi_0) =(\nu_0,1-\mu_0)$ holds for some $\xi_0\in M^{\rm c} $ with $\xi_0=(\nu_0,\mu_0)$.}

Then, by definition of $\phi_\sigma$, we have $M_{\xi_0}= M_{\xi_0}^{\rm c}=\emptyset$ and $(\xi_0, \xi_0)\not\in\rho$. 
Let $\xi=(\nu,\mu)\in M^{\rm c} $. Then, ${\rm sym}_\G(\rho)\supseteq \G(M)$ implies $((\nu,\mu),(\nu,1-\mu))\in\rho$ and,
by the consequences of (a), (b), and (c), $\{\xi\}\times (M_{\xi}\cup M_{\xi}^{\rm c}\cup\{\xi\mid (\xi, \xi)\in\rho\})=(\{\xi\}\times (I\setminus \{(\nu,1-\mu)\}))\cap \rho=\emptyset$ and thus, for all $\xi=(\nu,\mu)\in M^{\rm c}$, 
\[\phi_\sigma(\xi) =(\nu,1-\mu).\]
\end{list}
We can summarize the first four cases and characterize the choice function $\widetilde {\sigma_{M^{\rm c}}}\subseteq I^2$ for $\alpha \cap M^{\rm c}$ and $\rho\cap (M^{\rm c}\times I)$ given by $\phi_\sigma$ that is restricted to $M^{\rm c}$ as follows. Let $\xi_0$ be any individual in $M^{\rm c} $ and
\[\widetilde {\sigma_{M^{\rm c}}}=\left\{\begin{array}{ll} ((M^{\rm c}\setminus L)\times \{\zeta_0\} ) \cup ( L\times \{\zeta_1\})
&\mbox {if $\xi_0 $ satisfies (a),}\\

M^c \times \{\min_{\rm lexico}M_{\xi_0}\}
&\mbox {if $\xi_0 $ satisfies (b),}\\ 

 \{(\xi,\xi)\mid \xi\in M^c\}
&\mbox {if $\xi_0 $ satisfies (c),}\\ 

\{((\nu,\mu), (\nu,1-\mu))\mid \nu \in K^c\}
&\mbox {if $\xi_0 $ satisfies (d).}\end{array}\right .\] 
Second, let $M_\alpha=M\cap \alpha$ and let us assume that $M_\alpha$ is not empty. Then, let us consider the two disjoint subsets $M_1$ and $M_2$ given by $M_1=\{\xi\in M_ \alpha\mid \phi_\sigma(\xi) \in M^{\rm c}\}$ and $M_2=\{\xi\in M_ \alpha\mid \phi_\sigma(\xi) \in M\}$ and satisfying $M_\alpha=M_1\cup M_2$. Then, there are an $l=|M_\alpha|\geq 1$ and pairwise different individuals $\eta_1,\ldots,\eta_l\in M_\alpha$ satisfying $M_\alpha=\{\eta_1,\ldots,\eta_l\} $ and one of the three possible cases, $M_1=\{\eta_{1},\ldots,\eta_{t}\}$ and $M_2=\{\eta_{{t+1}},\ldots,\eta_{l}\}$ for some $t\in \{1,\ldots, l-1\}$, $M_1=\{\eta_1,\ldots,\eta_l\}$ and $M_2=\emptyset$, and $M_2=\{\eta_1,\ldots,\eta_l\}$ and $M_1=\emptyset$, respectively. Moreover, let $M=\{\eta_1,\ldots, \eta_k\}$ with $l\leq k=|M|=2|K|$ and $\eta_{s}\not= \eta_{r}$ for $1\leq s<r\leq k$.  For each single $\xi\in M_ \alpha$, only one of the following cases --- (e) and (f) --- is possible. 
\begin{list}{(\alph{li})}{\usecounter{li}}\setcounter{li}{4} 

\item {$\phi_\sigma(\xi) \in M^{\rm c}$ holds and thus $\xi\in M_1$.} Let $\xi=(\nu,\mu)$. Then, $K=K\cup \{\nu\}$ and, by Lemma \ref{ToIdealOfS2}, $\{\xi\}\times M^{\rm c}\subseteq \rho$ and thus $M^{\rm c}_{\xi}= M^{\rm c}$ and $\phi_\sigma(\xi)=\min_{\rm lexico}M^{\rm c}_{\xi}=\min{\!}_{\rm lexico}M^{\rm c}$ which means \[\phi_\sigma(\xi)=\zeta_0.\]
\item{$\phi_\sigma(\xi) \in M$ holds and thus $\xi\in M_2$.} 
In this case, we have, by definition, $\phi_\sigma(\xi)=\min_{\rm lexico}M_{\xi}$ or $\phi_\sigma(\xi)=\xi$ or $\phi_\sigma(\nu,\mu) =(\nu,1-\mu)$ for $ \xi=(\nu,\mu)$. Let $s=0$ if $M_1=\emptyset$ and otherwise $s=t$ and let $j_1, \ldots, j_{l-s}$ be the corresponding indices in $\{1, \ldots, k\}$ satisfying \[\phi_\sigma(\eta_{{s+1}})= \eta_{j_1}, \quad \ldots,\quad \phi_\sigma(\eta_{{l}})=\eta_{j_{l-s}}.\]
\end{list}
Let the choice functions $\widetilde{\sigma_{M_1}}\subseteq I^2$ and $\widetilde{\sigma_{M_2}}\subseteq I^2$ for $\alpha \cap M_1$ and $\rho\cap (M_1\times I)$ and for $\alpha \cap M_2$ and $\rho\cap (M_2\times I)$, respectively, be defined by
\[\widetilde{\sigma_{M_1}}=M_1\times \{\zeta_0\},\]
\[\widetilde{\sigma_{M_2}}= \left\{\begin{array}{ll} \{(\eta_{{t+1}},\eta_{j_1}), \ldots, (\eta_{{l}},\eta_{j_{l-t}})\} & \mbox{if $M_2\not=\emptyset$ and $M_1\not =\emptyset$}, \\ \{(\eta_{{1}},\eta_{j_1}), \ldots, (\eta_{{l}},\eta_{j_{l}})\} & \mbox{if $M_2\not=\emptyset$ and $M_1 =\emptyset$},\\\emptyset &\mbox{otherwise}.\end{array}\right.\]
For the predicate $\sigma$ defined by $\sigma (\xi,\eta) =true \Leftrightarrow \phi_\sigma(\xi)=\eta$ for all individuals $\xi,\eta\in I$ we have $\sigma= \sigma_{M^{\rm c}} \cup \sigma_{M_1} \cup \sigma_{M_2} $ and, because of ${\rm sym}_\G(\sigma_{M^{\rm c}})\supseteq \G(M\cup \{\zeta_0,\zeta_1\})$, ${\rm sym}_\G(\sigma_{M_1})\supseteq \G(M\cup \{\zeta_0\})$, and ${\rm sym}_\G(\sigma_{M_2})\supseteq \G(M)$, the inclusion ${\rm sym}_\G(\sigma)\supseteq \G(M\cup \{\zeta_0,\zeta_1\})$ and, thus, $\sigma\in J_2(I,\G,\I)$. \qed

\vspace{0.2cm}

 \begin{proposition} There is a formula $H$ in ${\cal L}^{(2)}_C$ such that 
 \[\S_2 \models \neg choice_*^1 (H).\]
 \end{proposition}

\noindent{\bf Proof.} We consider the domain $I$, the group $\G$, and the ideal $\I$ given in Definition \ref{Model2}, $\S_2=\S(I,\G,\I)$, $\S=\S_2$, and the symmetric and irreflexive relation $\tau=\{(\xi_1,\xi_2)\mid (\exists \nu\in \bbbn)(\{\xi_1,\xi_2\}=\{\nu \}\times N_1 )\}$. Because of ${\rm sym}_\G(\tau) \supseteq \G(\emptyset)$, $\tau$ belongs to $J_2(I,\G,\I)$. Let 
\[H(C)=\exists x\exists y (Cx\land Cy\land Txy\land \forall B (Bx\land By\land Txy \to \forall z(Cz\to Bz))).\]
Then, $\S\models_{\langle {C \atop\gamma} {T \atop\tau} \rangle}H(C)$ holds iff there is a $\nu\in \bbbn$ with $\gamma=\{\nu \}\times N_1 $. Consequently, we have $\S_{\langle{T\atop \tau}\rangle}(G_1\,\,\land \,\, G_2)=true$ for $ G_1=\forall C(H(C)\!\to\!\exists y Cy ) $ and $G_2= \forall C_1\forall C_2 (H(C_1)\land H(C_2)\land C_1\neq C_2 \to \neg \exists y (C_1y \land C_2y ))$.

In order to show that $choice_*^1(H)$ does not hold in $\S$, it is enough to prove that we have $\S_{\langle {T\atop \tau}\rangle}( \exists D\forall C(H(C)\to \exists !! y(Cy \land Dy )))=false$. For any $\delta\in J_1(I,\G,\I)$ there is a $K_\delta\subseteq \bbbn$ with ${\rm sym}_\G(\delta) \supseteq \G(K_\delta\times N_1 )$. For any $\nu\not \in K_\delta$, there are a $\pi \in \G(K_\delta\times N_1 )$, a $\phi_{\pi}\in \G_\bbbn$, and a $\pi_\nu\in N_1 $ such that $\phi_{\pi}(\nu)=\nu$, $\pi_\nu(0)=1$, $\pi_\nu(1)=0$, and thus $\pi((\nu,0))=(\phi_{\pi}(\nu),\pi_\nu(0))=(\nu,1)$ and $\pi((\nu,1))=(\phi_{\pi}(\nu),\pi_\nu(1))=(\nu,0)$ hold. This means that, 
we have $\delta \cap (\{\nu\}\times N_1 )=\emptyset$ or $\{\nu\}\times N_1 \subset \delta $ and, consequently, 
\[\S_{\langle{T\atop \tau}{C \atop\{\nu\}\times N_1 )} {D \atop\delta} \rangle}(H(C)\to \exists !! y(Cy \land Dy ))=false\] and thus 
\[\S_{\langle{T\atop \tau}{C \atop\{\nu\}\times N_1 )} {D \atop\delta} \rangle}(\neg(H(C)\to \exists !! y(Cy \land Dy )))=true.\]
This means, that we have 
\[\S_{\langle{T\atop \tau} \rangle}(\forall D\exists C(\neg(H(C)\to \exists !! y(Cy \land Dy )))=true,\]
\[\S_{\langle{T\atop \tau} \rangle}(\neg (\exists D\forall C(H(C)\to \exists !! y(Cy \land Dy )))=true,\] 
and
\[\S_{\langle{T\atop \tau} \rangle}(\exists D\forall C(H(C)\to \exists !! y(Cy \land Dy ))=false.$\qed$\] 

\vspace{0.2cm}
By the latter two propositions we get the following.
\begin{proposition}\label{ThACnotchoicedisj} There is a Henkin-Asser structure that is a model of $^{h}ax^{(2)} \cup \{AC^{1,1} \land \neg choice_*^1(H)\}$ for some $H$ in ${\cal L}^{(2)}_C$.
\end{proposition}

A similar statement is possible for $choice_h^1$ instead of $choice_*^1$ if we assume that ZFC is our metatheory (cf. Proposition \ref{TheoACnotChoice}). 

\subsection{Zermelo-Asser (1,1)-axiom does not imply Ackermann axioms}

 We assume that ZFC is consistent and that we have a model of ZFC (that we call {\em our metamodel}). For proving the next theorem, we want to consider a structure with an individual domain that is the union of two disjoint infinite sets. We combine a structure with the individual domain $\{0\}\times \bbbn$ and a structure with the individual domain $\{1\}\times \bbbn$ and permit all $n$-ary predicates $\alpha$ that are the finite union of $n$-ary predicates belonging to one of both structures (where we say that the predicate $\alpha$ is the union of the predicates $\alpha_1$ and $\alpha_2$ if $\widetilde \alpha=\widetilde\alpha_1\cup \widetilde\alpha_2$ holds). One of these substructures, the basic Henkin structure $\S_{\sf 0}^{\{0\}\times\bbbn}$, is isomorphic to $\S_0^\bbbn$ and, therefore, definable analogously to the basic Fraenkel structure $\S_0^\bbbn=\S(\bbbn,\G_{\sf 1} ^\bbbn,\I_{\sf 0})$. The second substructure, $\S_{\sf 1}^{\{1\}\times\bbbn}$, is the standard structure over $\{1\}\times \bbbn$ containing all $n$-ary predicates $\alpha$ with $\widetilde \alpha \subseteq (\{1\}\times \bbbn)^n$ (for any $n\geq 1$) in our metamodel. Whereas $\S_{\sf 0}^{\{0\}\times\bbbn}$ is a Henkin-Asser structure with a minimal set of predicates, $\S_{\sf 1}^{\{1\}\times\bbbn}$ is a Henkin-Asser structure with a maximal set of predicates within the framework of our metatheory.
 
 Let $I$ be the individual domain $\{0,1\}\times\bbbn$ and, moreover, let $\G=\{\pi\in \G_{\sf 1}^I\mid (\forall \xi\in\{1\}\times\bbbn)
(\pi(\xi)=\xi)\}$.

\begin{defi}[The model $\S_3$ of second order]\label{Model3} Let $\S_3$ be the Henkin-Asser structure $\S(I,\G,\I_0^{I})$.
\end{defi}

\begin{lemma}\label{Ideal3}$\G( K) =\G(K\cup P)=\G( K\cup (\{1\}\times \bbbn))$ holds for any finite set $K\subseteq \{0\}\times\bbbn$ and all $P\subseteq \{1\}\times \bbbn$.
\end{lemma}

\begin{proposition} $\S_3 \models AC^{1,1} $.\label{ACinS3}
\end{proposition}

\noindent{\bf Proof.} Let $\I=\I_0^I$ and $\S=\S_3$. Moreover, let $\rho$ and $\alpha$ be predicates in $J_2(I,\G,\I)$ and $J_1(I,\G,\I)$, respectively, such that $\S\models _{\langle{R\atop \rho}{A\atop\alpha}\rangle} \forall x (Ax \leftrightarrow \exists y Rx y)$ is satisfied. For defining $\sigma$, we construct a choice function $\phi_\sigma$ for $\rho$ and $\alpha$ in a similar way as in the proof of Proposition \ref{ACinS2}. Let $K $ be a finite set with $K\subseteq \{0\}\times\bbbn$ and ${\rm sym}_\G(\rho)\supseteq\G( K )$ which exists by Lemma \ref{Ideal3}. Let $L=(\{1\}\times \bbbn)\cup K$, $L^{\rm c}=I\setminus L $ which implies $L^{\rm c}\subseteq \{0\}\times\bbbn$, $\eta_0=\min_{\rm lexico} L^{\rm c}$, $M=L^{\rm c} \setminus \{\eta_0\}$, and $\eta_1=\min_{\rm lexico}M$.
For all $\xi\in \alpha$, let 
\begin{eqnarray*} L_{\xi}\hspace{0.1cm}&=&\{\eta\in\hspace{0.1cm}L\hspace{0.13cm}\mid (\xi,\eta )\in \rho\},\\
N_ {\xi}\hspace{0.09cm}&=&\{\eta \in L^{\rm c} \mid (\xi,\eta )\in \rho\,\,\&\,\,\xi\not=\eta \}, \end{eqnarray*}
and
\[\phi(\xi)=\left\{\begin{array}{ll} \min_{\rm lexico} L_{\xi} & \mbox{ if } L_{\xi}\not= \emptyset,\\\min_{\rm lexico} N_ {\xi} & \mbox{ if } L_{\xi}= \emptyset \mbox{ and } N_ {\xi}\not= \emptyset, \\ 
 \xi& \mbox{ if } L_{\xi}= \emptyset \mbox{ and } N_ {\xi}= \emptyset\mbox{ and } (\xi, \xi)\in\rho.\end{array}\right.\]By this, the choice function $\phi:\widetilde\alpha\to I$ defines a choice relation $\sigma$ for $\alpha$ and $\rho$. Let $\sigma$ be the predicate in ${\rm pred}_2(I)$ given by $\widetilde \sigma=\phi$ and let $\phi_\sigma=\phi$. The domain $\widetilde\alpha$ of $\phi_\sigma$ is the union $\alpha_1\cup\cdots \cup \alpha_5$ given by the five predicates $ \alpha_1 =\{\xi \in L^{\rm c}\mid \phi_\sigma(\xi) \in L_\xi\}$, $ \alpha_2=\{\xi \in L\mid \phi_\sigma(\xi) \in L_\xi\}$, $ \alpha_3=\{\xi \in L^{\rm c}\mid \phi_\sigma(\xi) \in N_ {\xi}\}$, $ \alpha_4=\{\xi \in L\mid \phi_\sigma(\xi) \in N_ {\xi}\}$, and $ \alpha_5=\{\xi \in L^{\rm c}\mid \phi_\sigma(\xi) =\xi\}$. If $\phi_\sigma(\xi) =\xi$ and $\xi\in L$, then $\xi$ is also in $ L_\xi$. Hence, we need only consider the following five cases. For showing that $\sigma$ belongs to $\S_3$, let us prove for every $i\in\{1,\ldots, 5\}$ that the predicate $\sigma_i$ given by $\sigma_i=\sigma \cap (\alpha_i \times I)$ is a choice relation for $\alpha_i$ and $\rho_i=\rho \cap (\alpha_i \times I)$ that belongs to $\S_3$. 

\begin{list}{(\alph{li})}{\usecounter{li}}
\item {$\phi_\sigma(\xi_0) \in L_{\xi_0}$ for some $\xi_0\in I$.} 

Then, $\alpha_1\cup \alpha_2\not=\emptyset$. 

\newcounter{iii}
\begin{list}{(\alph{li}\arabic{iii})}{\usecounter{iii}}
\item Let $\alpha_1\not=\emptyset$ and $\xi_0\in \alpha_1$ be fixed. Moreover, let $\xi\in L^{\rm c}$. Then, for each $\eta\in L$ there is a $\pi\in \G(K)$ with $\pi(\eta)=\eta$, $\pi(\xi)=\xi_0$, and $\pi(\xi_0)=\xi$. Therefore, we have $L_{\xi}=L_{\xi_0}$ and thus $\xi\in \alpha_1$. Let $L_0= L_{\xi_0}$ and $\eta_2=\min_{\rm lexico} L_0$ which means, for instance, $\phi_\sigma(\xi_0)=\eta_2$. Consequently, $\phi_\sigma(\xi)= \min_{\rm lexico} L_{\xi}=\eta_2$ holds for all $\xi\in L^{\rm c}$ and the choice relation $\sigma_1$ for $\alpha_1$ and $\rho_1$ is given by
\[\sigma_1=\{(\xi, \eta_2)\mid \xi\in L^{\rm c}\,\,\&\,\, \phi_\sigma(\xi)=\eta_2\}\]
 and we have ${\rm sym}_{\G}(\sigma_1 )\supseteq \G(K\})$.

\item Let $\alpha_2\not=\emptyset$ and $\xi\in \alpha_2$. Then, $\xi\in L$ and $(\pi(\xi), \pi(\phi_\sigma(\xi))=(\xi,\phi_\sigma(\xi))$ for all $\pi\in \G(K)$. From this observation we can conclude that the choice relation $\sigma_2$ for $\alpha_2$ and $\rho_2$ is given by
\[\sigma_2=\{(\xi, \phi_\sigma(\xi))\mid \xi\in \alpha_2 \,\,\&\,\, \phi_\sigma(\xi)\in L \}\]
and by $\alpha_2\subseteq L$ we have we have ${\rm sym}_{\G}(\sigma_2)\supseteq \G(K)$.
\end{list}

\item{$\phi_\sigma(\xi_0) \in N_ {\xi_0}$ for some $\xi_0\in I$.} 

Then, $\alpha_3\cup \alpha_4\not=\emptyset$. Let $\xi\in \alpha_3\cup \alpha_4$. Then, by definition, $\phi_\sigma(\xi)\not=\xi$ and, for each $\eta\in L^{\rm c}\setminus \{\xi\}$, there is a $\pi\in \G(K)$ with $\pi(\phi_\sigma(\xi))=\eta$ and $\pi(\xi)=\xi$ which means that $\{\xi\}\times (L^{\rm c}\setminus \{\xi\}) \subseteq\rho$ holds and, thus, $\phi_\sigma(\xi)=\eta_0$ holds if $\xi\not = \eta_0$ and otherwise $\phi_\sigma(\xi)=\eta_1$ holds. 

\begin{list}{(\alph{li}\arabic{iii})}{\usecounter{iii}}

\item Let $\xi_0\in \alpha_3$. Then, $\xi_0\in L^{\rm c}$. Moreover, first, let $\xi$ be any individual in $M $. If $\xi_0\in M$, then there is a $\pi\in \G(K)$ with $\pi(\eta_0)=\eta_0$ and $\pi(\xi_0)=\xi$. If $\xi_0=\eta_0$, then there is a $\pi\in \G(K)$ with $\pi(\eta_1)=\eta_0$ and $\pi(\eta_0)=\xi$. Thus, for $\xi\in M$, we have $(\xi,\eta_0)\in \rho$ and $\phi_\sigma(\xi)=\eta_0$ and thus $\xi\in \alpha_3$. If $\xi=\eta_0$ and $\xi_0\in M$, then there is a $\pi\in \G(K)$ with $\pi(\eta_0)=\eta_1$ and $\pi(\xi_0)=\eta_0$. This means that $(\eta_0,\eta_1)\in \rho$ and $\phi_\sigma(\eta_0)=\eta_1$ and thus $\eta_0\in \alpha_3$. Consequently, if $\alpha_3\not=\emptyset$, then $\alpha_3=L^c$. Therefore, the choice relation $\sigma_3$ for $\alpha_3$ and $\rho_3$ is given by 
\[\sigma_3=\{(\xi, \eta_0)\mid \xi\in M \,\,\&\,\, \phi_\sigma(\xi_0)=\eta_0\}\cup \{(\eta_0,\eta_1)\mid \phi_\sigma(\eta_0)=\eta_1\}\]
and we have ${\rm sym}_{\G}(\sigma_3)\supseteq \G(K\cup \{\eta_0,\eta_1\})$.

\item If $\xi\in \alpha_4$, then $\xi\in L$. Thus, 
\[\sigma_4 =\{(\xi, \eta_0)\mid \xi\in \alpha_4\,\,\&\,\, \phi_\sigma(\xi)= \eta_0\}\]
 is a choice relation for $\alpha_4$ and $\rho_4$ and by $\alpha_4\subseteq L$ we have ${\rm sym}_{\G}(\sigma_4)\supseteq \G(K\cup\{\eta_0\})$.
\end{list}

\item {$\phi_\sigma(\xi_0) =\xi_0$ for some $\xi_0\in I$.} 

Then, $L_{\xi_0}=N_ {\xi_0}=\emptyset$ and, thus, $\xi_0\in L^{\rm c}$. This means that $(\xi, \xi)\in\rho$ holds for all $\xi\in L^{\rm c}$ and, thus, $\phi_\sigma(\xi) =\xi$ holds for all $\xi\in L^{\rm c}$ with $L_{\xi}= N_ {\xi}=\emptyset$.

\begin{list}{(\alph{li}\arabic{iii})}{\usecounter{iii}}
\item $\xi_0\in L^{\rm c}$ means that $\alpha_5\not=\emptyset$. Consequently,
\[\sigma_5 =\{(\xi, \xi)\mid \xi\in L^{\rm c} \,\,\&\,\, (\{\xi\}\times \{\eta_0,\eta_1,\eta_2\})\cap\rho=\emptyset \}\]
 is a choice relation for $\alpha_5$ and $\rho_5$ and we have ${\rm sym}_{\G}(\sigma_5 )\supseteq \G(K\cup \{\eta_0,\eta_1\})$.

\item For $\xi_0\in L$ see case (a2).
\end{list}
\end{list}

\noindent
Consequently, we have ${\rm sym}_\G(\sigma)\supseteq \G((K \cup \{\eta_0,\eta_1\})$ for $\sigma=\sigma_1\cup \cdots \cup \sigma_5 $ and, thus, $\sigma\in J_2(I,\G,\I)$ and $\sigma$ is a choice relation for $\alpha$ and $\rho$.\qed

\vspace{0.3cm}

As in \cite[p.\,37]{1c}, we describe the property of a predicate to be a bijection by a formula. Let the formula $bij( R,A,D)$\label{DefBij} be defined by
\[bij( R,A,D)=\forall x \forall y (Rxy\to Ax\land Dy)\land \forall y(Dy \to\exists!!x Rxy)\]\[\hfill \land\,\forall x(Ax \to\exists!!y Rxy).\hspace*{1.9cm}\] Thus, for any predicate structure $\S$, $\Sigma\models_{\langle{A\atop\alpha } {D\atop\delta } { R\atop \rho} \rangle}bij( R,A,D)$ means that $\widetilde\rho$ is (the graph of) a bijective function $\phi: \widetilde\alpha \to \widetilde \delta$. In the following, to simplify matters, we say that such a predicate $\rho$ is a {\em bijection with the domain $\alpha$ and the codomain $\delta$}. We will use that each predicate $\rho$ for which $\widetilde\rho$ is a function $\phi: \widetilde\alpha \to \widetilde \delta$ on a finite domain $\widetilde\alpha$ belongs to each Henkin-Asser structure whose domain of individuals includes $\widetilde\alpha $ and $\widetilde \delta$. The image $\widetilde\gamma_\xi$ of $\{\xi\}$ under any relation $\widetilde\rho$, given by $\gamma_\xi(\eta)= \Sigma_{\langle { R\atop \rho} {x \atop\xi}{y\atop\eta } \rangle} (Rxy)$, will be described by the formula $image(x,R,C)$ given by 
\[image(x,R,C)= \forall y (Rxy\leftrightarrow Cy).\]
 Moreover, since we consider again ZFC as metatheory, there are a well-ordering on the whole set $\bbbn$, which we denote by $\leq$, and elements denoted by $0,1,\ldots$ such that $0\leq 1 \leq 2\leq \cdots $ holds.

\begin{proposition}
There exists a formula $H$ in ${\cal L}^{(2)}_{x,D}$ such that we have 
\[\S_3 \models \neg choice_h^{1,1}(H).\]
\end{proposition}
\noindent{\bf Proof.} We consider the domain $I$, the group $\G$, and the ideal $\I$ given in Definition \ref{Model3} such that $\S_3=\S(I,\G,\I)$. Let $\S=\S_3$, $\alpha_0=\{0\}\times \bbbn$, $\alpha_1=\{1\}\times \bbbn$, $\tau =\{((1,\mu_1),(1,\mu_2))\in \alpha_1 ^2\mid \mu_2\leq\mu_1\}$, $f$ be any assignment in $\S$, and $f'=f\langle{A_0\atop\alpha_0 } {A_1\atop\alpha_1 } {T\atop \tau} \rangle$. We will show that $\S_{f'}(choice_h^{1,1}(H))$ is $false$ for the Ackermann axiom $choice_h^{1,1}(H)$ given by $ \forall x \exists D H(x,D) \to \exists S \forall x \exists D (\forall y(D y \leftrightarrow Sxy) \land H(x, D ))$ if $H(x,D)$ is the formula $\exists R G(x,D,R)$ given by 
\[ G(x,D,R)=A_1x \to \forall y(Dy\to A_0y) \land \forall C( image(x,T,C) \to bij(R,C,D)). \]
 For every $n\in \bbbn$, let $\xi_n=(1,n) $ and $\rho_n$ be a binary predicate given by $\rho_n=\{((1,\mu),(0,\mu))\mid \mu\leq n\}\subseteq \alpha_1 \times\alpha_0$.
Thus, $\rho_n$ is a bijection with a finite domain and $\{(1,\mu)\in \alpha_1\mid \mu \leq n\}$ is the domain of $\rho_n$ as well as the image of $\{\xi_n\}$ under $\tau$. For every $n\in \bbbn$, let $\delta_n$ be a unary predicate given by $\widetilde\delta_n=\{(0,0),(0,1), \ldots, (0,n)\}$. Then, $ \delta_n$ is the codomain of $\rho_n$. Consequently, for any $n\in\bbbn$, $\Sigma\models_{f'\langle {x\atop \xi_n}{D\atop \delta_{n} } {R\atop \rho_{n}} \rangle}G(x,D,R)$ holds. We can say that, for any $\xi\in \alpha_1$, there are $\rho$ and $\delta$ such that we have $\Sigma\models_{f'\langle{x\atop \xi} {D\atop \delta} {R\atop \rho} \rangle}G(x,D,R)$. This means that we get 
\[\S\models_{f'}\forall x \exists D H(x,D)\] 
and --- in any case --- for each predicate $\delta$ satisfying $\S_{f'\langle{D\atop \delta}{x\atop \xi_{n}}\rangle} ( H(x,D))=true$ for some $n\geq 0$, the condition $|\delta|=n+1$ holds.

Now, let us assume that $\S$ is a model of $ choice_h^{1,1}$. Then, by this assumption, there is a $\sigma\in J_2(I,\G,\I)$ such that
\[\S\models_{f'\langle\,{S\atop\sigma}\rangle}\forall x \exists D (\forall y(D y \leftrightarrow Sxy) \land H(x,D)).\] 
 The assumption $\sigma \in J_2(I,\G,\I)$ implies by Lemma \ref{Ideal3} that there is a finite subset $K_\sigma \subseteq \bbbn$ such that ${\rm sym}_\G(\sigma)\supseteq\G(\{0\} \times K_\sigma)$. Let $\nu=|K_\sigma|$ and let $\delta$ be the image of $\{\xi_\nu\}$ under $\widetilde \sigma$ which means $\delta=\sigma_{\xi_\nu}=_{\rm df}\{\eta \mid (\xi_\nu,\eta)\in \sigma\}$. This $\delta$ is the only predicate satisfying $\S\models_{f'\langle{D\atop \delta}{x\atop \xi_{\nu}}{S\atop\sigma}\rangle} \forall y(D y \leftrightarrow Sxy)$. Thus, $\S\models_{f'\langle{D\atop \delta}{x\atop \xi_{\nu}}\rangle} H(x,D)$ must also be true for this $\delta$.
Because of the latter property we have $\delta \subseteq \alpha_0$ and $|\delta|= \nu+1$ and thus $|(\{\xi_\nu\}\times \alpha_0)\cap \sigma|=\nu+1$. The latter fact means that there is a $\mu_0\in \bbbn\setminus K_\sigma$ such that $(\xi_\nu,(0,\mu_0))\in\sigma$. Let $\mu$ be any number in $ \bbbn\setminus K_\sigma$. Since there is a $\pi\in {\rm sym}_\G(\sigma)$ with $\pi((0,\mu_0))=(0,\mu)$ and $\pi(\xi_\nu)=\xi_\nu$, the relationships $(\xi_\nu,(0,\mu))\in\sigma$ and thus $(0,\mu)\in \delta$ for the corresponding $\delta=\sigma_{\xi_\nu}$ hold. However, then $|\delta|>\nu+1$ and this contradicts $|\delta|=\nu+1$. Therefore, our assumption is wrong and, thus, $\S\models \neg choice_h^{1,1}(H)$. 

\qed 
\begin{corollary} 
There is a formula $H$ in ${\cal L}^{(2)}_{x,D}$ such that 
$\S_3 \models \neg choice^{1,1}(H)$.
\end{corollary}
 
Provided that ZFC is consistent, we get the following theorem.

\begin{proposition}\label{TheoACnotChoice} There is a Henkin-Asser structure that is a model of $^{h}ax^{(2)} \cup \{AC^{1,1} \land \neg choice_h^{1,1}(H)\}$ for some $H$ in ${\cal L}^{(2)}_{x,D}$.
\end{proposition}

\section {A summary}\label{Summary}

\setcounter{satz}{0}

In \,{\rm HPL}, $AC_*^{1,1}$ is weaker than $AC^{n,m}$ for any $n,m\geq 1$. $\S_0$ is a model of $AC^{1,1}$ (cf.\,\,\cite{Gass94}) and a model of $choice_h^{1,1}$ (cf.\,\,\cite{Gass24B}). Moreover, it can be shown that $\S_0$ is a model of $choice_*^1$. Thus, the main results of this paper are the following theorems. 

\begin{theorem}\label{HPL_ind_1} $choice_*^1$ is {\rm HPL}-independent of $AC^{1,1}$. 
\end{theorem}

We know that there are a formula $H$ in ${\cal L}^{(2)}$ and a Henkin-Asser structure that is a model of $AC^{1,1}\land \neg choice_*^1(H)$. $\S_0$ is a model of $choice_*^1\cup \{AC^{1,1}\}$.

\begin{theorem} \label{HPL_ind_2} $choice_h^{1,1}$ is {\rm HPL}-independent of $AC^{1,1}$. 
\end{theorem}

We know that there are a formula $H$ in ${\cal L}^{(2)}$ and a Henkin-Asser structure that is a model of $AC^{1,1}\land \neg choice_h^{1,1}(H)$ if ZFC is consistent. $\S_0$ is a model of $choice_h^{1,1}\cup \{AC^{1,1}\}$.

\newpage

\begin{figure}[ht] \unitlength1cm

\begin{picture}(12,7.9) \thicklines
\put(0,0.5){
\put(-3,2.0)
{
\put(5.75,3.38){\vector(1,0){1.5}}
\put(7.25,3.02){\vector(-1,0){1.5}}
\put(5.75,0.68){\vector(1,0){1.5}}
\put(7.25,0.32){\vector(-1,0){1.5}}
\put(4.95,3.2){\oval(1.3,1)\makebox(0,0){$AC^{n,m}$}}
\put(4.95,0.5){\oval(1.3,1)\makebox(0,0){$AC^{1,m}$}}
}
\put(1.75,-0.7)
{
\put(4.25,6.08){\vector(1,0){1.5}}
\put(5.75,5.72){\vector(-1,0){1.5}}
\put(4.25,3.38){\vector(1,0){1.5}}
\put(5.75,3.02){\vector(-1,0){1.5}}
\put(6.32,5.3){\vector(0,-1){1.5}}
\put(6.68,3.8){\vector(0,1){1.5}}
\put(6.45,4.5){\line(1,0){0.48}}

\put(7.1,5.35){\vector(1,-1){1.75}}
\put(9.1,3.85){\vector(-1,1){1.75}}
\put(8.2,4.5){\line(0,1){0.48}}

\put(7.1,2.65){\vector(1,-1){1.75}}
\put(9.1,1.15){\vector(-1,1){1.75}}
\put(8.2,1.8){\line(0,1){0.48}}

\put(9.5,2.6){\vector(0,-1){1.5}}
\put(8.75,3.2){\vector(-1,0){1.5}}
\put(8.22,2.95){\line(0,1){0.48}}

\put(3.4,5.9){\oval(1.5,1)\makebox(0,0){$AC^{n,n+m}_{*}$}}
\put(3.4,3.2){\oval(1.5,1)\makebox(0,0){$AC^{1,n+m}_{*}$}}
\put(6.5,5.9){\oval(1.4,1)\makebox(0,0){$AC^{n,n+1}_{*}$}}
\put(6.5,3.2){\oval(1.3,1)\makebox(0,0){$AC^{1,2}_{*}$}}
\put(9.5,3.2){\oval(1.3,1)\makebox(0,0){$AC^{n,1}_{*}$}}
\put(9.5,0.5){\oval(1.3,1)\makebox(0,0){$AC^{1,1}_{*}$}}

\put(-0.7,7.6){\oval(2,1)\makebox(0,0){$choice_*^{n+m}$}}
\put(-1,7){\vector(0,-1){1.6}}
\put(-0.1,7){\vector(1,-2){0.25}}

\put(-0.7,4.8){\oval(2,1)\makebox(0,0){$choice_*^{2}$}}
\put(-1,4.2){\vector(0,-1){3}}
\put(-0.1,4.2){\vector(1,-2){0.23}}

\put(-0.7,0.65){\oval(2,1)\makebox(0,0){$choice_*^{1}$}}
\put(0.3,2.6){\vector(-1,-2){0.7}}
\put(-0.2,2){\line(1,0){0.45}}

\put(0.45,0.65){\vector(1,0){8.3}}
}
}
\end{picture}

\vspace{1cm}

\begin{picture}(12,7.9) \thicklines
\put(0,0.5){
\put(-3,2.0)
{
\put(5.75,3.38){\vector(1,0){1.5}}
\put(7.25,3.02){\vector(-1,0){1.5}}
\put(5.75,0.68){\vector(1,0){1.5}}
\put(7.25,0.32){\vector(-1,0){1.5}}

\put(4.95,3.2){\oval(1.3,1)\makebox(0,0){$AC^{n,m}$}}
\put(4.95,0.5){\oval(1.3,1)\makebox(0,0){$AC^{1,m}$}}
}
\put(1.75,-0.7)
{
\put(4.25,6.08){\vector(1,0){1.5}}
\put(5.75,5.72){\vector(-1,0){1.5}}
\put(4.25,3.38){\vector(1,0){1.5}}
\put(5.75,3.02){\vector(-1,0){1.5}}
\put(6.32,5.3){\vector(0,-1){1.5}}
\put(6.68,3.8){\vector(0,1){1.5}}
\put(6.45,4.5){\line(1,0){0.48}}
\put(7.1,5.35){\vector(1,-1){1.75}}
\put(9.1,3.85){\vector(-1,1){1.75}}
\put(8.2,4.5){\line(0,1){0.48}}
\put(7.1,2.65){\vector(1,-1){1.75}}
\put(9.1,1.15){\vector(-1,1){1.75}}
\put(8.2,1.8){\line(0,1){0.48}}
\put(9.5,2.6){\vector(0,-1){1.5}}
\put(8.75,3.2){\vector(-1,0){1.5}}
\put(8.22,2.95){\line(0,1){0.48}}
\put(3.4,5.9){\oval(1.5,1)\makebox(0,0){$AC^{n,n+m}_{*}$}}
\put(3.4,3.2){\oval(1.5,1)\makebox(0,0){$AC^{1,n+m}_{*}$}}
\put(6.5,5.9){\oval(1.4,1)\makebox(0,0){$AC^{n,n+1}_{*}$}}
\put(6.5,3.2){\oval(1.3,1)\makebox(0,0){$AC^{1,2}_{*}$}}
\put(9.5,3.2){\oval(1.3,1)\makebox(0,0){$AC^{n,1}_{*}$}}
\put(9.5,0.5){\oval(1.3,1)\makebox(0,0){$AC^{1,1}_{*}$}}
\put(-0.7,7.6){\oval(2,1)\makebox(0,0){$choice_h^{n,m}$}}
\put(-1,7){\vector(0,-1){1.6}}
\put(-0.1,7){\vector(1,-2){0.25}}
\put(-0.7,4.8){\oval(2,1)\makebox(0,0){$choice_h^{1,1}$}}
\put(-0.1,4.2){\vector(1,-2){0.23}}
\put(-0.58,3.9){\line(1,0){0.4}}
\put(-0.3,3.7){\vector(-1,2){0.28}}
}
}
\end{picture}
\caption{Outlook: Possible extensions of Fig. 1 ($n>1$, $m\geq 1$, incomplete)}
\end{figure}
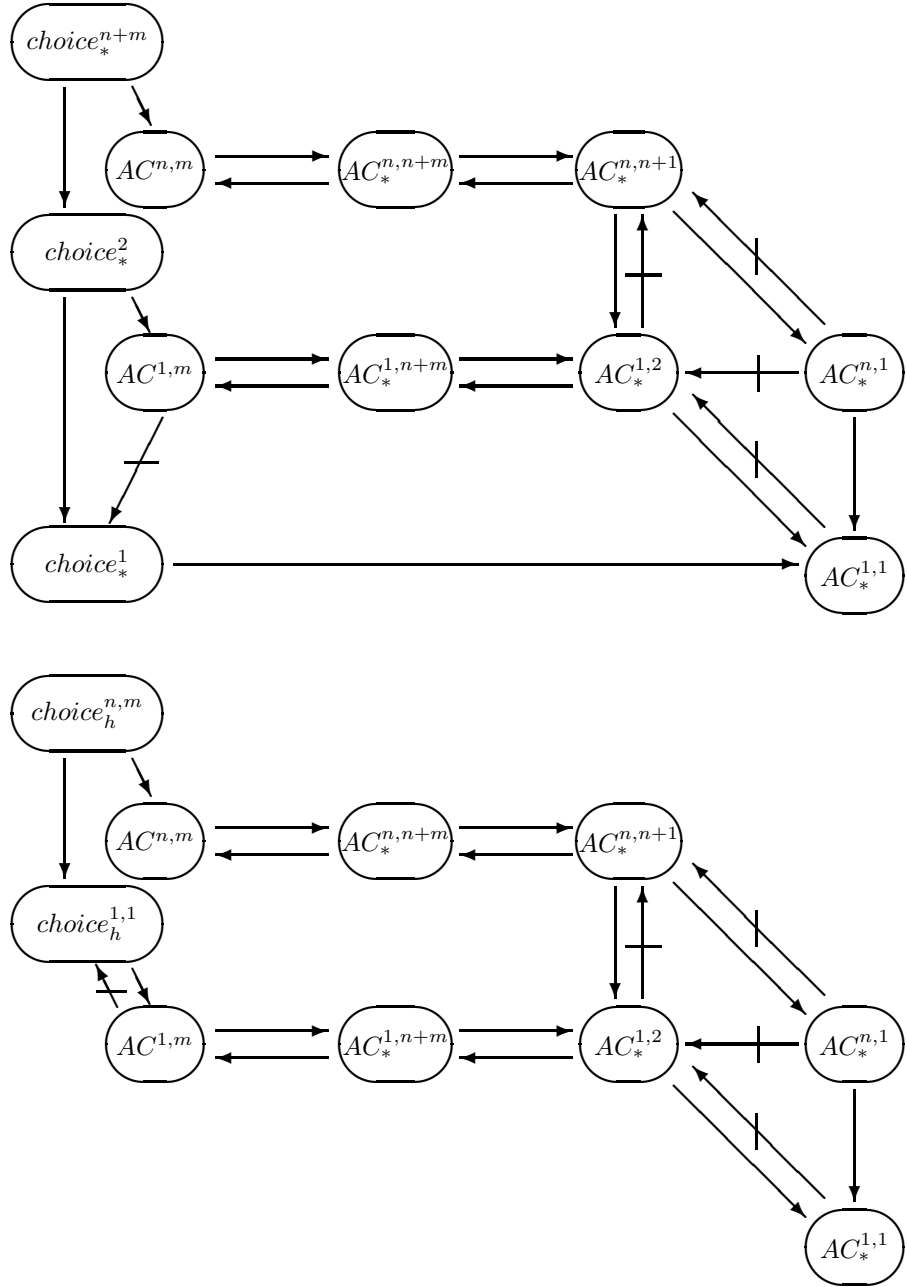
\end{document}